\documentclass[12pt]{article}
\usepackage{amsmath,amssymb,amsthm,comment}

\newtheorem{theorem}{Theorem}[section]
\newtheorem{thmy}{Theorem}

\newtheorem{corollary}[theorem]{Corollary}

\newtheorem{example}[theorem]{Example}

\newcommand{\dd}{\displaystyle }

\def\barr{\begin{array}}
\def\earr{\end{array}}

\title{The Reidemeister spectrum of ZM-groups}
\author{Marius T\u arn\u auceanu}
\date{November 13, 2025}

\begin{document}

\maketitle

\begin{abstract}
Given a group $G$ and an automorphism $\varphi$ of $G$, two elements $x,y\in G$ are said to be $\varphi$-conjugate if $x=gy\varphi(g)^{-1}$ for some $g\in G$. The number $R(\varphi)$ of equivalence classes with respect to this relation is called the Reidemeister number of $\varphi$ and the set $\{R(\varphi)|\varphi\in {\rm Aut}(G)\}$ is called the Reidemeister spectrum of $G$. In this paper, we determine the Reidemeister spectrum of ZM-groups, extending some results of \cite{13}.
\end{abstract}

{\small
\noindent
{\bf MSC 2020\,:} Primary 20D45, 20E45; Secondary 20E22. 

\noindent
{\bf Key words\,:} ZM-group, twisted conjugacy, Reidemeister number, Reidemeister spectrum.} 

\section{Introduction}

Let $G$ be a group and ${\rm Aut}(G)$ be the group of automorphisms of $G$. The concept of conjugacy in $G$ can be naturally generalized to the so-called
\textit{twisted conjugacy}: given $\varphi\in {\rm Aut}(G)$, two elements $x,y\in G$ are called \textit{$\varphi$-conjugate} if there exists $g\in G$ such that $x=gy\varphi(g)^{-1}$. We define the \textit{Reidemeister number} $R(\varphi)$ of $\varphi$ as the number of equivalence classes with respect to this relation. Note that $R(\varphi)\in\mathbb{N}_0\cup\{\infty\}$. Also, we define the \textit{Reidemeister spectrum} of $G$ to be ${\rm Spec}_R(G)=\{R(\varphi)|\varphi\in {\rm Aut}(G)\}$. If ${\rm Spec}_R(G)=\{\infty\}$, we say that $G$ has the $R_{\infty}$-property.

Twisted conjugacy arises naturally in Nielsen fixed-point theory (see e.g. \cite{12}), but there is also a strong algebraic interest in its study, in particular in the study of $R_{\infty}$-property (see e.g. \cite{5,8,9,15}). For many groups either partial or full information is known about their Reidemeister spectrum. It has either been proven that they have the $R_{\infty}$-property (such as Baumslag-Solitar groups \cite{7} or Thompson’s group \cite{1}), or that they do not (such as certain groups of exponential growth \cite{10} or free groups of infinite rank \cite{3}), and for some even the complete Reidemeister spectrum has been determined (such as low-dimensional crystallographic groups \cite{4} or finite abelian groups \cite{14}).

The starting point for our discussion is given by the paper \cite{13}, where the Reidemeister spectrum has been determined for split metacyclic groups of the form $(C_n\times C_{p^m})\rtimes C_p$ where $p$ is prime, $m$ and $n$ are non-negative integers and $n$ is coprime with $p$. In what follows we will deal with ZM-groups, i.e. finite groups with all Sylow subgroups cyclic. By \cite{11} (Ch. IV, Satz 2.11.), such a group is of type
\begin{equation}
{\rm ZM}(m,n,r)=\langle a, b \mid a^m = b^n = 1,
\hspace{1mm}b^{-1} a b = a^r\rangle, \nonumber
\end{equation}
where the triple $(m,n,r)$ satisfies the conditions
\begin{equation}
(m,n)=(m,r-1)=1\,\,\footnote{\,\,For a finite number of positive integers $a_1$, $a_2$, \dots, $a_n$, we denote by $(a_1,a_2,\dots,a_n)$ their greatest common divisor.} \mbox{ and } r^n
\equiv 1 \hspace{1mm}({\rm mod}\hspace{1mm}m). \nonumber
\end{equation}It is clear that $|{\rm ZM}(m,n,r)|=mn$ and $Z({\rm ZM}(m,n,r))=\langle b^{o_m(r)}\rangle$, where
\begin{equation}
o_m(r)={\rm min}\{k\in\mathbb{N}^* \mid r^k\equiv 1 \hspace{1mm}({\rm mod} \hspace{1mm}m)\}\,\,\footnote{\,\,For simplicity of notation, throughout the paper we will use $d$ instead of $o_m(r)$.}\nonumber
\end{equation}is the multiplicative order of $r$ modulo $m$. For every non-negative integer $u$, we define
\begin{equation}
[u]_r=\left\{\barr{lll}
    \!\!1+r+\cdots+r^{u-1},&u>0\\
    \!\!0,&u=0.\earr\right.\nonumber
\end{equation}

Our main result is stated as follows.

\begin{theorem}
We have
\begin{equation}
{\rm Spec}_R({\rm ZM}(m,n,r))=\left\{\!\dd\frac{1}{m}\dd\sum_{v=0}^{m-1}\sum_{^{\,0\leq u\leq n-1}_{\frac{n}{(n,y-1)}\mid u}}\frac{(m,[u]_r)}{o_{\frac{(m,[u]_r)}{(m,[u]_r,v)}}(r)}\,\mid\,0\leq y<n, y\equiv 1\, ({\rm mod}\, d)\!\right\}\!.\nonumber
\end{equation}
\end{theorem}

\begin{example}
For the dicyclic group ${\rm Dic}_3={\rm ZM}(3,4,2)={\rm SmallGroup}(12,1)$, we have $d=2$ and the above formula leads to
\begin{equation}
{\rm Spec}_R({\rm Dic}_3)=\{4,6\}.\nonumber
\end{equation}
\end{example}\newpage

A particular case of Theorem 1.1 is the following.

\begin{corollary}
If $n$ is prime, then we have
\begin{equation}
{\rm Spec}_R({\rm ZM}(m,n,r))=\left\{n-1+\dd\frac{S}{n}\right\}\!,\nonumber
\end{equation}where
\begin{equation}
S=\dd\sum_{u=0}^{n-1}\,(m,[u]_r).\nonumber
\end{equation}
\end{corollary}

\begin{example}
For $n=2$ and $r=m-1$ we have ${\rm ZM}(m,n,r)\cong D_{2m}$, the dihedral group of order $2m$. Then $S=m+1$ and so the formula in Corollary 1.3 becomes
\begin{equation}
{\rm Spec}_R(D_{2m})=\left\{\dd\frac{m+3}{2}\right\}, \mbox{ for all odd integers } m\geq 3.\nonumber
\end{equation}
\end{example}

We remark that a split metacyclic group $(C_n\times C_{p^m})\rtimes C_p$ becomes a ZM-group for $m=0$ and so Theorem 1.1 extends the results obtained in \cite{13} in this particular case.

Our approach uses some basic tools of group action theory, such as the Burnside's lemma, and involves modular arithmetic.

\bigskip\noindent{\bf Burnside's lemma.}
{\it Let $G$ be a finite group acting on a finite set $X$ and
\begin{equation}
Fix(g)=\{x\in X \mid g\circ x=x\}, \mbox{ for all } g\in G.\nonumber
\end{equation}Then the number of distinct orbits is
\begin{equation}
k=\frac{1}{|G|}\sum_{g\in G}|Fix(g)|.\nonumber
\end{equation}}
 
Most of our notation is standard and will not be repeated here. Basic definitions and results on group theory can be found in \cite{11}.

\section{Proofs of the main results}

The general form of automorphisms of a metacyclic group has been determined in the main result of \cite{2}. For ZM-groups, in \cite{6} we obtained that:

\begin{theorem}
Each automorphism of ${\rm ZM}(m,n,r)$ is given by
\begin{equation}
b^ua^v\mapsto b^{yu}a^{x_1v+x_2[u]_r},\, u,v\geq 0,\nonumber
\end{equation}for a unique triple of integers $(x_1,x_2,y)$ such that
\begin{equation}
0\leq x_1,x_2<m,\, (x_1,m)=1,\, 0\leq y<n \mbox{ and } y\equiv 1 \hspace{1mm}({\rm mod}\hspace{1mm}d).
\end{equation}In particular, we have
\begin{equation}
|{\rm Aut}({\rm ZM}(m,n,r))|=m\varphi(m)\frac{n}{d}\,.\nonumber
\end{equation}
\end{theorem}

In what follows we will denote by $f_{(x_1,x_2,y)}$ the automorphism of ${\rm ZM}(m,n,r)$ determined by a triple $(x_1,x_2,y)$ satisfying (1).

We are now able to prove our results.

\begin{proof}[Proof of Theorem 1.1.] 
Let $f=f_{(x_1,x_2,y)}\in{\rm Aut}({\rm ZM}(m,n,r))$. In order to compute $R(f)$, we will apply Burnside's lemma to the action of ${\rm ZM}(m,n,r)$ on ${\rm ZM}(m,n,r)$ given by
\begin{equation}
g\circ x=gxf(g)^{-1}, \forall\, g,x\in{\rm ZM}(m,n,r).\nonumber
\end{equation}Let $g=b^ua^v$ and $x=b^{\alpha}a^{\beta}$, where $0\leq u,\alpha<n$ and $0\leq v,\beta<m$. Using Theorem 2.1, we get
\begin{align*}
Fix(g)&=\{x\in{\rm ZM}(m,n,r)\mid f(g)=x^{-1}gx\}\\
&=\{b^{\alpha}a^{\beta}\in{\rm ZM}(m,n,r)\mid b^{yu}a^{x_1v+x_2[u]_r}=b^ua^{r^{\alpha}v-\beta(r^u-1)}\}\\
&=\{b^{\alpha}a^{\beta}\in{\rm ZM}(m,n,r)\mid b^{yu}=b^u \mbox{ and } a^{x_1v+x_2[u]_r}=a^{r^{\alpha}v-\beta(r^u-1)}\}.\nonumber
\end{align*}We have
\begin{equation}
b^{yu}=b^u\Leftrightarrow n\mid (y-1)u\Leftrightarrow\frac{n}{(n,y-1)}\mid u\nonumber
\end{equation}and
\begin{equation}
a^{x_1v+x_2[u]_r}=a^{r^{\alpha}v-\beta(r^u-1)}\Leftrightarrow m\mid (r^{\alpha}-x_1)v-\beta(r^u-1)-x_2[u]_r.\nonumber
\end{equation}For fixed $u,v$ and $\alpha$, the last congruence has solutions $\beta$ if and only if 
\begin{equation}
(m,[u]_r)\mid (r^{\alpha}-x_1)v.
\end{equation}Moreover, if it has solutions, then there are exactly $(m,[u]_r)$ solutions. We also remark that (2) is equivalent with
\begin{equation}
\frac{(m,[u]_r)}{(m,[u]_r,v)}\mid r^{\alpha}-x_1
\end{equation}and (3) has the same number of solutions $\alpha$ as
\begin{equation}
\frac{(m,[u]_r)}{(m,[u]_r,v)}\mid r^{\alpha}-1,\nonumber
\end{equation}that is as
\begin{equation}
o_{\frac{(m,[u]_r)}{(m,[u]_r,v)}}(r)\mid\alpha.
\end{equation}Obviously, (4) has $\frac{n}{o_{\frac{(m,[u]_r)}{(m,[u]_r,v)}}(r)}$ solutions and so we get
\begin{align*}
\hspace{5mm}R(f)&=\dd\frac{1}{mn}\dd\sum_{g\in{\rm ZM}(m,n,r)}|Fix(g)|\\
&=\dd\frac{1}{mn}\dd\sum_{(u,v)}|\{(\alpha,\beta)\mid b^{yu}=b^u \mbox{ and } a^{x_1v+x_2[u]_r}=a^{r^{\alpha}v-\beta(r^u-1)}\}|\\
&=\dd\frac{1}{mn}\sum_{v=0}^{m-1}\sum_{^{\,0\leq u\leq n-1}_{\frac{n}{(n,y-1)}\mid u}}\sum_{\alpha=0}^{n-1}|\{\beta=\overline{0,m-1}\,\mid\, m\!\mid\! (r^{\alpha}-x_1)v-\beta(r^u-1)-x_2[u]_r\}|\\
&=\dd\frac{1}{mn}\sum_{v=0}^{m-1}\sum_{^{\,0\leq u\leq n-1}_{\frac{n}{(n,y-1)}\mid u}}\!\sum_{^{\,\,\,\,\,\,\,\,\,\,\,0\leq \alpha\leq n-1}_{\,\,\,(m,[u]_r)\mid (r^{\alpha}-x_1)v}}\!\!\!(m,[u]_r)\\
&=\dd\frac{1}{mn}\sum_{v=0}^{m-1}\sum_{^{\,0\leq u\leq n-1}_{\frac{n}{(n,y-1)}\mid u}}(m,[u]_r)|\{\alpha=\overline{0,n-1}\,\mid\, (m,[u]_r)\!\mid\! (r^{\alpha}-x_1)v\}|\\
&=\dd\frac{1}{mn}\sum_{v=0}^{m-1}\sum_{^{\,0\leq u\leq n-1}_{\frac{n}{(n,y-1)}\mid u}}(m,[u]_r)\,\frac{n}{o_{\frac{(m,[u]_r)}{(m,[u]_r,v)}}(r)}\\
&=\dd\frac{1}{m}\dd\sum_{v=0}^{m-1}\sum_{^{\,0\leq u\leq n-1}_{\frac{n}{(n,y-1)}\mid u}}\frac{(m,[u]_r)}{o_{\frac{(m,[u]_r)}{(m,[u]_r,v)}}(r)}\,,\nonumber
\end{align*}as desired.
\end{proof}

\begin{proof}[Proof of Corollary 1.2.] 
If $n$ is prime, then we have $d=n$ and so the conditions $0\leq y<n$ and $y\equiv 1\, ({\rm mod}\, d)$ imply $y=1$ for any $f=f_{(x_1,x_2,y)}\in{\rm Aut}({\rm ZM}(m,n,r))$. Thus, we get
\begin{equation}
R(f)=\dd\frac{1}{m}\dd\sum_{v=0}^{m-1}\sum_{u=0}^{n-1}\frac{(m,[u]_r)}{o_{\frac{(m,[u]_r)}{(m,[u]_r,v)}}(r)}\,.\nonumber
\end{equation}Since $o_{\frac{(m,[u]_r)}{(m,[u]_r,v)}}(r)$ is a divisor of $n$, we have either $o_{\frac{(m,[u]_r)}{(m,[u]_r,v)}}(r)=1$ or $o_{\frac{(m,[u]_r)}{(m,[u]_r,v)}}(r)=n$. We remark that
\begin{equation}
o_{\frac{(m,[u]_r)}{(m,[u]_r,v)}}(r)=1\Leftrightarrow\frac{(m,[u]_r)}{(m,[u]_r,v)}=1\Leftrightarrow (m,[u]_r)\mid v.\nonumber
\end{equation}It follows that
\begin{align*}
R(f)&=\dd\frac{1}{m}\dd\sum_{v=0}^{m-1}\sum_{u=0}^{n-1}\frac{(m,[u]_r)}{o_{\frac{(m,[u]_r)}{(m,[u]_r,v)}}(r)}\\
&=\dd\frac{1}{m}\dd\sum_{u=0}^{n-1}(m,[u]_r)\left(\sum_{^{\,\,0\leq v\leq m-1}_{\,(m,[u]_r)\,\mid\, v}}1+\sum_{^{\,\,0\leq v\leq m-1}_{\,(m,[u]_r)\,\nmid\, v}}\frac{1}{n}\right)\\
&=\dd\frac{1}{m}\dd\sum_{u=0}^{n-1}(m,[u]_r)\left[\frac{m}{(m,[u]_r)}+\frac{1}{n}\left(m-\frac{m}{(m,[u]_r)}\right)\right]\\
&=\dd\sum_{u=0}^{n-1}\left(1-\frac{1}{n}\right)+\frac{1}{n}\sum_{u=0}^{n-1}\,(m,[u]_r)\\
&=n-1+\dd\frac{S}{n}\,,\nonumber
\end{align*}completing the proof.
\end{proof}

\bigskip\noindent{\bf Acknowledgements.} The author is grateful to the reviewer for remarks which improve the previous version of the paper.

\vspace*{3ex}\small

\hfill
\begin{minipage}[t]{5cm}
Marius T\u arn\u auceanu \\
Faculty of  Mathematics \\
``Al.I. Cuza'' University \\
Ia\c si, Romania \\
e-mail: {\tt tarnauc@uaic.ro}
\end{minipage}

\end{document}